\documentstyle{amsart}

\title{Dense Egyptian Fractions}
\author{Greg Martin}
\address{School of Mathematics\\Institute for Advanced Study\\Olden
Lane\\Princeton, NJ 08540\\U.S.A.}
\email{gerg@@math.ias.edu}

\newtheorem{theorem}{Theorem}
\newtheorem{lemma}[theorem]{Lemma}

\newenvironment{pflike}[1]{\noindent{\bf #1}}{\vskip10pt} 
\newenvironment{proof}{\begin{pflike}{Proof:}}{\qed\end{pflike}}

\newcommand{\abs}[1]{\left|#1\right|}
\newcommand{\floor}[1]{\left\lfloor#1\right\rfloor}

\renewcommand{\mod}[1]{{\ifmmode\text{\rm\ (mod~$#1$)}\else\discretionary{}{}{\hbox{ }}\rm(mod~$#1$)\fi}}
\newcommand{\ep}{\varepsilon}
\renewcommand{\implies}{\Rightarrow}

\newcommand{\half}{{\mathchoice{\textstyle\frac12}{1/2}{1/2}{1/2}}}
\newsymbol\dnd 232D

\renewcommand{\lg}[1]{\mathop{\log_{#1}}}
\def\lgs#1^#2{\mathop{\log_{#1}^{#2}}}

\newcommand{\A}{{\cal A}}
\newcommand{\B}{{\cal B}}
\newcommand{\C}{{\cal C}}
\newcommand{\D}{{\cal D}}

\renewcommand{\S}{{\cal S}}
\newcommand{\T}{{\cal T}}
\newcommand{\Z}{{\cal Z}}

\newcommand{\lol}[2]{{\mathchoice{\log#1\over\log#2}{\log#1/\log#2}{\log#1/\log#2}{\log#1/\log#2}}}
\newcommand{\err}{{\frac{\lgs2^3x}{\log y}}}

\newcommand{\abfrac}[1]{{\mathchoice{a_{#1}\over b_{#1}}{a_{#1}/b_{#1}}{}{}}}
\newcommand{\abcfrac}[2]{{\mathchoice{a_{#1,#2}\over b_{#1,#2}}{a_{#1,#2}/b_{#1,#2}}{}{}}}
\newcommand{\cdfrac}[2]{{\mathchoice{a'_{#1,#2}\over b'_{#1,#2}}{a'_{#1,#2}/b'_{#1,#2}}{}{}}}

\begin{document}

\abstract{Every positive rational number has representations as
Egyptian fractions (sums of reciprocals of distinct positive integers)
with arbitrarily many terms and with arbitrarily large denominators.
However, such representations normally use a very sparse subset of the
positive integers up to the largest denominator. We show that for
every positive rational there exist representations as Egyptian
fractions whose largest denominator is at most $N$ and whose
denominators form a positive proportion of the integers up to $N$, for
sufficiently large $N$; furthermore, the proportion is within a small
factor of best possible.}\endabstract

\maketitle

\section{Introduction}

The ancient Egyptians wrote positive rational numbers as sums of
distinct reciprocals of positive integers, or unit fractions. In 1202,
Fibonacci published an algorithm (subsequently rediscovered by
Sylvester in 1880, among others) for constructing such
representations, which have come to be called Egyptian fractions, for
any positive rational number. Since that time, number theorists have
been interested in some quantitative aspects of Egyptian fraction
representations. For instance, there are algorithms that improve upon
the Fibonacci--Sylvester algorithm in various ways, bounding the size
of the largest denominator or limiting the number of terms. Bleicher
\cite{Ble:ANAftEoEF} has a thorough survey of, and references to, such
developments.

One line of questions concentrates on the number of terms in Egyptian
fraction representations. A positive rational $m/n$ can always be
expanded into an Egyptian fraction with at most $m$ terms, for
instance by the Farey series algorithm (see Golomb
\cite{Gol:AAAftRPotAP}). Erd\H{o}s and Straus conjectured that for all
integers $n\ge2$, the fraction $4/n$ can actually be written as the
sum of three unit fractions rather than four. (In their conjecture,
distinctness of the terms is not required.) Sierpi\'nski
\cite{Sie:SlDdNReFP} made the same conjecture for fractions of the
form $5/n$, and mentioned that Schinzel conjectured that for any fixed
numerator $m$, the fraction $m/n$ could be written as the sum of three
unit fractions for sufficiently large $n$. In this vein, Vaughan
\cite{Vau:OaPoESaS} showed that almost all positive fractions with
numerator $m$ can be written as the sum of three unit fractions.

One can also investigate the behavior of the number of terms in
Egyptian fractions at the other extreme. Any unit fraction can be
split into two using the identity $1/n=1/(n+1)+1/(n(n+1))$;
consequently, given a particular Egyptian fraction representation, one
can repeatedly use this splitting algorithm on the term with largest
denominator to construct such representations with arbitrarily many
terms. However, this process results in a tremendously thin set of
denominators in the sense that if the largest denominator is $x$,
then the number of denominators is $\ll\log\log x$. One can try to use
the splitting algorithm on the intermediate terms, but then the danger
arises that the resulting unit fractions will no longer be distinct.

It is therefore interesting to ask how many terms can be used to
represent a rational number as an Egyptian fraction, given a bound on
the size of the largest denominator. The purpose of this paper is to
demonstrate that a positive proportion of the integers up to the bound
can in fact be assembled to form such a representation. We will establish
the following theorem:

\begin{theorem}
Let $r$ be a positive rational number and $\eta>0$ a real number. For
any real number $x$ that is sufficiently large in terms of $r$ and
$\eta$, there is a set $\S$ of positive integers not exceeding $x$
such that $r = \sum_{n\in\S} 1/n$ and $\abs\S > (C(r)-\eta)x$, where
\begin{equation*}
C(r) = (1-\log2)\Big({ 1-\exp\Big( {-r\over1-\log2} \Big) }\Big).
\end{equation*}
\label{densethm}
\end{theorem}

That one can always find such a dense set of denominators is already
surprising (though the Egyptians might not have been pleased to do
their arithmetic this way!), but it turns out that the proportion
given in Theorem \ref{densethm} is comparable to the best possible
result. Given a positive rational number $r$, any set $\S$ of positive
integers not exceeding $x$ with cardinality $\floor{(1-e^{-r})x}$
satisfies
\begin{equation}
\sum_{n\in\S} 1/n \ge\!\! \sum_{e^{-r}x+1<n\le x} \!\!1/n = r +
O_r(x^{-1}).
\label{derive}
\end{equation}
Consequently, the best possible value for $C(r)$ in Theorem
\ref{densethm} would be $1-e^{-r}$. Notice that
$C(r)/(1-e^{-r})=1-O(r)$ as $r$ tends to zero, and that
\begin{equation*}
{C(r)\over1-e^{-r}} = (1-\log2){1-\exp(-r/(1-\log2))\over1-e^{-r}} >
1-\log 2 = 0.30685\dots.
\end{equation*}
Thus when $r$ is small, the value for $C(r)$ given in Theorem
\ref{densethm} is very nearly best possible, and in any case it is
never smaller than 30\% of the best possible value.

While the proof of Theorem \ref{densethm} involves a rather complex
notation, the idea underlying the construction is quite
straightforward. One begins by subtracting from $r$ the reciprocals of
a suitable set $\A$ of integers not exceeding $x$ with cardinality at
least $(C(r)-\eta)x$. For a given prime $p$ dividing the denominator
of this difference, we can add back in the reciprocals of a very few
members of $\A$ so that the factors of $p$ in the denominator are
cancelled out. If we repeat this process for all large primes, we will
have expressed $r$ as the sum of the reciprocals of almost all of the
elements of $\A$, plus a small rational number whose denominator is
only divisible by small primes. We are then able to write this small
rational number as the sum of reciprocals of distinct integers much
smaller than the members of $\A$, using the idea outlined above
combined with an existing algorithm for expanding rational numbers
into Egyptian fractions.

For a given prime $p$ and rational number $r$ whose denominator is a
multiple of $p$, we might have to add the reciprocals of as many as
$p-1$ multiples of $p$ to eliminate the factor of $p$ from the
denominator of $r$, as shown in Lemmas \ref{Zpsums} and
\ref{addfracs}. Since all of the elements of $\A$ are at most $x$ in
size, we must have $p(p-1)\le x$; thus we must restrict our attention
to integers that are roughly $x^{1/2}$-smooth.  Fortunately, the
distribution of smooth numbers has been widely studied. For instance,
the density of $x^{1/2}$-smooth integers of size $x$ is $1-\log2$ (see
Section \ref{sns}); this is the origin of such factors in the
expression for $C(r)$.

Moreover, it turns out that we must consider integers that are
divisible by different powers of $p$ separately, which forces us to
have at least $p-1$ multiples of $p^2$ in our set $\A$, at least $p-1$
multiples of $p^3$, etc. Thus we must restrict to $x^{1/2}$-smooth
integers that are squarefree with respect to primes exceeding
$x^{1/3}$, cubefree with respect to primes exceeding $x^{1/4}$, and so
on. We shall find that the extra conditions on the multiplicity of
their prime factors is easily handled. For simplicity we shall work
with, roughly, $x^{1/2}$-smooth numbers that are $k$-free (not
divisible by the $k$th power of any prime) and that are squarefree
with respect to primes exceeding $x^{1/k}$, for some integer
$k\ge2$. Such integers satisfy all the above constraints, and we show
in Lemma \ref{smoothsets} that the multiples of each prime power
among such integers are sufficiently plentiful for our argument to
go through.

We define $\lg1x=\max\{\log x,1\}$ and $\lg jx=\lg1\,(\lg{j-1}x)$ for
any integer $j\ge2$, and we write $\log^kx$ and $\lgs j^kx$ as
shorthand for $(\log x)^k$ and $(\lg jx)^k$ respectively. We use the
notations $P(n)$ and $p(n)$ to denote the greatest and least prime
factors of $n$ respectively, and make the conventions that $P(1)=1$
and $p(1)=\infty$. We say that a prime power $p^l$ exactly divides an
integer $n$, or that $n$ is exactly divisible by $p^l$, if $p^l$
divides $n$ but $p^{l+1}$ does not. The constants implicit in the
$\ll$ and $O$-notations in this paper may depend on any Greek variable
($\alpha$, $\delta$, $\ep$, $\eta$, and $\lambda$) and also on $k$ and
the rational number $r$ or $a/b$ where appropriate, but they will
never depend on $p$, $q$, $t$, $v$, $w$, $x$, or $y$; this remains the
case when any of these variables is adorned with primes or
subscripts. When the phrase ``$x$ is sufficiently large'' is used, the
size of $x$ may depend on the Greek variables and $k$, $r$, and $a/b$
as well.

The author would like to thank Hugh Montgomery and Trevor Wooley for
their guidance and support and the referee for many valuable
comments. This material is based upon work supported under a National
Science Foundation Graduate Research Fellowship.

\section{Elementary Lemmas}\label{elemsec}

The following lemma is an easy consequence of the Cauchy--Davenport
Theorem (see \cite[Lemma 2.14]{Vau:THLM}, for instance), but we
provide a direct proof.

\begin{lemma}
Let $t$ be a nonnegative integer, and let $x_1, \dots, x_t$ be
nonzero elements of ${\bf Z}_p$, not necessarily distinct. Then the
number of elements of ${\bf Z}_p$ that can be written as the sum of
some subset (possibly empty) of the $x_i$ is at least
$\min\{p,t+1\}$. In particular, if $t\ge p-1$, then every element of
${\bf Z}_p$ can be so written.
\label{Zpsums}
\end{lemma}

\noindent We remark that the conclusion of the lemma is best possible,
since we could take $x_1\equiv\dots\equiv x_t\not\equiv0\mod p$, in
which case the set of elements of ${\bf Z}_p$ that can be written as
the sum of a subset of the $x_i$ is precisely $\{0,x,\dots,tx\}$,
which has cardinality $t+1$ if $t<p$.

\begin{proof}
We use induction on $t$, the case $t=0$ being trivial. Given a
positive integer $t$ and nonzero elements $x_1, \dots, x_t$ of
${\bf Z}_p$, let $\S_i$ be the set of elements of ${\bf Z}_p$ that can
be written as the sum of some subset of $x_1, \dots, x_i$; then
certainly $\S_i\subset\S_t$ and so $\abs{\S_i}\le\abs{\S_t}$. By
induction we may assume that $\abs{\S_{t-1}}\ge\min\{p,t\}$. If
$\abs{\S_{t-1}}=p$ then $\abs{\S_t}=p$ as well, so we may assume that
$p>\abs{\S_{t-1}}\ge t$.

Suppose that $\abs{\S_t}<t+1$. Then we have $\abs{\S_t}\le t\le
\abs{\S_{t-1}}\le\abs{\S_t}$, and so
$\abs{\S_t}=\abs{\S_{t-1}}=t<p$. In particular, the map
$f:\S_{t-1}\to\S_t$ defined by $f(y)=y+x_t$ is actually a bijection of
$\S_{t-1}$ onto itself. Thus
\begin{equation*}
\sum_{y\in\S_{t-1}}y \equiv \sum_{y\in\S_{t-1}}(y+x_t) \equiv tx_t +
\sum_{y\in\S_{t-1}}y \mod p,
\end{equation*}
which implies that $tx_t \equiv 0\mod p$, a contradiction since $x_t$
is nonzero\mod p and $0<t<p$. Therefore $\abs{\S_t}\ge t+1$, which
establishes the lemma.
\end{proof}

Using this lemma, we can show that a factor of a prime $p$ can be
eliminated from the denominator of a rational number by adding the
reciprocals of fewer than $p$ integers from any prescribed set meeting
certain conditions.

\begin{lemma}
Let $p^l$ be a prime power, and let $N$ be an integer that is exactly
divisible by $p^l$. Let $c/d$ be a rational number with $d$ dividing
$N$, and let $\S$ be a set of integers dividing $N$ that are all
exactly divisible by $p^l$. If $\abs\S\ge p-1$, then there is a subset
$\T$ of $\S$ with cardinality less than $p$ such that, if we define
\begin{equation*}
\frac{c'}{d'} = \frac cd + \sum_{n\in\T} \frac1n
\end{equation*}
with $c'/d'$ in lowest terms, then $d'$ divides $N/p$.
\label{addfracs}
\end{lemma}

\begin{proof}
Without loss of generality we may assume that $\abs\S=p-1$, whereupon
we denote the elements of $\S$ by $n_1, \dots, n_{p-1}$. Let
$M=\mathop{\rm lcm}\{d, n_1, \dots, n_{p-1}\}$, let $m=M/d$, and let
$m_i=M/n_i$ for each $1\le i\le p-1$. Note that $M$ divides $N$ and
that each $m_i$ is nonzero\mod p. By Lemma \ref{Zpsums}, every element
of ${\bf Z}_p$ can be written as the sum of some subset of $m_1,
\dots, m_{p-1}$. In particular, we can choose distinct indices $i_1$,
\dots, $i_t$ for some $0\le t\le p-1$ so that
\begin{equation}
m_{i_1}+\dots+m_{i_t}\equiv -cm\mod p.
\label{pickedright}
\end{equation}
If we set $\T=\{n_{i_1}, \dots, n_{i_t}\}$, then
\begin{equation*}
\frac cd + \sum_{n\in\T} \frac1n = \frac{cm}M + \sum_{j=1}^t
\frac{m_{i_j}}M = {(cm+m_{i_1}+\dots+m_{i_t})/p\over M/p},
\end{equation*}
where the numerator is an integer by virtue of equation
(\ref{pickedright}). Since $M$ divides $N$, this establishes the
lemma.
\end{proof}

We now cite an existing Egyptian fraction algorithm, which we will use
near the end of the proof of Theorem \ref{densethm}.

\begin{lemma}
Let $c/d$ be a positive rational number with $d$ odd. Suppose that
$c/d<1/P(d)$. Then there exists a set $\C$ of distinct odd positive
integers, with
\begin{equation*}
\max\{n\in\C\} \ll d\!\!\prod_{p\le P(d)}\!\!p,
\end{equation*}
such that $c/d=\sum_{n\in\C}1/n$.
\label{Breusch}
\end{lemma}

\begin{proof}
Breusch has shown that any positive rational number with odd
denominator can be written as the sum of reciprocals of distinct odd
positive integers. When one examines his construction \cite[Lemmas
1--3]{Breusch}, one finds that when $c/d<1/P(d)$, the integers
involved in fact do not exceed
\begin{equation*}
5\mathop{\rm lcm}\{d,\,3^2\!\!\!\!\prod_{3<p\le P(d)}\!\!\!\!p\},
\end{equation*}
which implies the bound given in the statement of the lemma.
\end{proof}

\section{Smooth Number Sets}\label{sns}

For real numbers $x$, $y>1$, we recall the definition of $\A(x,y)$,
the set of $y$-smooth numbers up to $x$:
\begin{equation*}
\A(x,y) = \{n\le x:P(n)\le y\}.
\end{equation*}
Let $w>1$ and $0\le\lambda<1$ be real numbers and $k\ge2$ be an
integer. We will need to work with the following sets of smooth
numbers with various specified properties:
\begin{equation}
\begin{split}
\A(x,y;w,\lambda) &= \{n: \lambda x<n\le x;\, P(n)\le y;\, n\hbox{ is
$k$-free};\, d^2\mid n\implies P(d)\le w\}, \\
\A(x,y;w,\lambda;p^l) &= \{n\in\A(x,y;w,\lambda): n=mp^l\hbox{ with }
P(m)<p\}, \\
\A_0(x,y;w,\lambda) &= \{n\in\A(x,y;w,\lambda): n\hbox{ is odd, not
of the form }m^2+m-1\}, \\
\A_0(x,y;w,\lambda;p^l) &= \A_0(x,y;w,\lambda) \cap
\A(x,y;w,\lambda;p^l). \\
\end{split}
\label{smsetsdef}
\end{equation}
The first of these sets is the set of $k$-free smooth numbers that are
squarefree with regard to large primes, as described in the
introduction, while the second is the subset of the first consisting
of those integers whose largest prime factor is $p$ with multiplicity
exactly $l$. The sets $\A(x,y;w,\lambda;p^l)$ and
$\A_0(x,y;w,\lambda;p^l)$ do not actually depend on the parameter $y$
as long as $y\ge p$; but we retain the parameter $y$ for consistency
of notation, and as a reminder that
$\A(x,y;w,\lambda;p^l)\subset\A(x,y;w,\lambda)$ and likewise for the
$\A_0$-sets.

We note that for any $y'<y$, we may write $\A(x,y;w,\lambda)$ as the
disjoint union
\begin{equation*}
\A(x,y;w,\lambda) = \A(x,y';w,\lambda) \cup \bigg( \bigcup
\begin{Sb}y'<p\le y \\ p>w\end{Sb} \A(x,y;w,\lambda;p) \bigg)
\cup \bigg( \bigcup\begin{Sb}y'<p\le y \\ p\le w\end{Sb}
\bigcup_{l=1}^{k-1} \A(x,y;w,\lambda;p^l) \bigg),
\end{equation*}
and the same is true if we replace every occurrence of $\A$ by $\A_0$.

As is customary, we let $\Psi(x,y)$ denote the cardinality of
$\A(x,y)$. We will also use $\Psi$, with any list of arguments and
with or without subscript, to denote the cardinality of the
corresponding $\A$-set. It is well-known that $\Psi(x,y)\sim\rho(\log
x/\log y)x$ for a certain range of $x$ and $y$, where the Dickman
rho-function $\rho(u)$ is defined for $u>0$ to be the continuous
solution to the differential-difference equation
\begin{equation}
\begin{split}
\rho(u)=1, \quad& 0<u\le1; \\
u\rho'(u)=-\rho(u-1), \quad& u>1. \\
\end{split}
\label{rhodef}
\end{equation}
The following lemma, due to Hildebrand, describes this asymptotic
formula more specifically. We have not made any effort to optimize the
error term or the range of $y$ in the hypothesis, as it will suffice
for our purposes as stated.

\begin{lemma}
Let $x\ge 3$ and $y$ be real numbers satisfying $\exp(\lgs2^4x)\le
y\le x$. Then
\begin{equation*}
\Psi(x,y) = x\rho \Big( \lol xy \Big) \Big({ 1 + O\Big( {\lg2x\over\log
y} \Big) }\Big).
\end{equation*}\vskip-20pt
\label{hild}
\end{lemma}

\vskip12pt\begin{proof}
Hildebrand shows \cite[Theorem 1]{Hil:smooth} that for any $\ep>0$,
one has
\begin{equation*}
\Psi(x,x^{1/u}) = x\rho(u) \Big({ 1 + O\Big( {u\lg1u\over\log x} \Big) }\Big)
\end{equation*}
uniformly for $x\ge3$ and $1\le u\le\log x/(\lg2x)^{5/3+\ep}$. On
setting $\ep=7/3$ and $u=\log x/\log y$ (so that $y=x^{1/u}$), the
lemma follows immediately.
\end{proof}

From the definition (\ref{rhodef}) of $\rho(u)$, one can easily derive
that $\rho$ is a positive, decreasing function, and that
$\rho(u)=1-\log u$ for $1\le u\le2$, so that $\rho(2)=1-\log2$ in
particular. We will need the following additional properties of
$\rho(u)$.

\begin{lemma}
Let $x\ge3$, $v\ge1$, and $y$ be real numbers satisfying $v\le x$ and
$\exp(\lgs2^4x)\le y\le x$. Then:
\begin{enumerate}
\item[\rm(a)]{
$\displaystyle\rho \Big( \lol{x/v}y \Big) \ll \rho \Big( \lol xy \Big)
v^{\lg2x/\log y}$;
}
\item[\rm(b)]{
for any $\ep>0$, we have $\displaystyle\rho \Big( \lol xy \Big)^{\!\!-1}
\ll x^\ep$;
}
\item[\rm(c)]{
for any real number $\alpha>1$, we have
\begin{equation*}
\sum_{\log x<n\le x^{1/\alpha}} \!\!n^{-\alpha} \rho \Big(
\lol{xn^{-\alpha}}y \Big) \ll \rho \Big( \lol xy \Big) (\log
x)^{-\alpha+1};
\end{equation*}
}
\item[\rm(d)]{
if $\displaystyle\log v=o\bigg({\log y\over\lgs2^2x}\bigg)$, then
$\displaystyle\rho \Big( \lol{x/v}y \Big) = \rho \Big( \lol xy \Big)
\Big({ 1 + O\Big( {\log v\lgs2^2x\over\log y} \Big) }\Big)$.  }
\end{enumerate}
\label{rholemma}
\end{lemma}

\begin{proof}
Hildebrand shows \cite[Lemma 1(vi)]{Hil:smooth} that for any $0\le
t\le u$, we have
\begin{equation}
\rho(u-t) \ll \rho(u) (u\lgs1^2u)^t.  \label{hildd}
\end{equation}
Part (a) follows immediately on taking $u=\lol xy$ and $t=\lol vy$
and noting that $u\lgs1^2u\le\log x$. Part (b) follows from part (a)
on taking $v=x$ and noting that
\begin{equation*}
x^{\lg2x/\log y} \le x^{1/\lgs2^3x} \ll x^\ep.
\end{equation*}

By using part (a) again, the sum in part (c) is
\begin{equation}
\ll \rho \Big( \lol xy \Big) \!\sum_{\log x<n\le x^{1/\alpha}}
\!n^{-\alpha} (n^\alpha)^{\lg2x/\log y} \ll \rho \Big( \lol xy \Big)
\!\sum_{n>\log x} \!n^{-\alpha+\alpha/\lgs2^3x}.
\label{oddnpower}
\end{equation}
For $x$ sufficiently large, the exponent $-\alpha+\alpha/\lgs2^3x$ is
bounded above away from $-1$. Therefore the right hand side of
(\ref{oddnpower}) is
\begin{equation*}
\begin{split}
&\ll \rho \Big( \lol xy \Big) (\log x)^{-\alpha+\alpha/\lgs2^3x+1} \\
&= \rho \Big( \lol xy \Big) (\log x)^{-\alpha+1} \exp (
\alpha/\lgs2^2x ) \ll \rho \Big( \lol xy \Big) (\log
x)^{-\alpha+1},
\end{split}
\end{equation*}
which establishes part (c). 

Finally, for any real numbers $1<s<t$, we have
\begin{equation*}
\begin{split}
\rho(s) - \rho(t) &= -\int^t_s \rho'(u)\,du = \int^t_s
\frac{\rho(u-1)}u\,du \\
&\ll \frac{\rho(s-1)}s(t-s) \ll \rho(s)(s\lgs1^2s)^1 \Big(\frac{t-s}s\Big)
\end{split}
\end{equation*}
by equation (\ref{hildd}). Therefore
\begin{equation*}
\rho(t) = \rho(s)(1+O((t-s)\lgs1^2s)).
\end{equation*}
Letting $s=\lol{(x/v)}y$ and $t=\lol xy$, we see that
\begin{equation*}
\rho \Big( \lol xy \Big) = \rho \Big( \lol{x/v}y \Big) \Big({ 1 +
O\Big( {\log v\lgs2^2x\over\log y} \Big) }\Big);
\end{equation*}
and under the hypothesis that $\log v=o(\log y/\lgs2^2x)$, the
$(1+O(\cdot))$-term can be inverted and moved to the other side of the
equation, which establishes part (d).

\end{proof}

With Lemmas \ref{hild} and \ref{rholemma} at our disposal, we can
establish the following lemmas concerning the distributions of the
sets of smooth numbers defined in equation (\ref{smsetsdef}).

\begin{lemma}
Let $x\ge 3$, $y$, and $w$ be real numbers satisfying
$\exp(\lgs2^4x)\le y\le x$ and $w\ge\log x$. Then:
\begin{enumerate}
\item[\rm(a)]{
$\displaystyle \Psi(x,y;w,0) = \frac x{\zeta(k)} \rho \Big( \lol xy \Big)
\Big({ 1 + O\Big( \err \Big) }\Big)$;
}
\item[\rm(b)]{
$\displaystyle \Psi_0(x,y;w,0) = \frac x{\xi(k)} \rho \Big( \lol xy \Big)
\Big({ 1 + O\Big( \err \Big) }\Big)$, where $\xi(k)=2(1-2^{-k})\zeta(k)$.
}
\end{enumerate}
\label{kfrees}
\end{lemma}

\vskip12pt\begin{proof}
From the definition (\ref{smsetsdef}) of $\A(x,y;w,\lambda)$, we
have
\begin{equation*}
x^{-1}\Psi(x,y;w,0) = x^{-1} \!\!\sum \begin{Sb}n\le x \\ P(n)\le
y\end{Sb} \!\!\bigg(\! \sum \begin{Sb}d^k\mid n \\ P(d)\le w\end{Sb}
\mu(d) \bigg) \bigg(\!  \sum \begin{Sb}f^2\mid n \\ p(f)>w\end{Sb}
\mu(f) \bigg).
\end{equation*}
Interchanging the order of summation yields
\begin{equation*}
\begin{split}
x^{-1}\Psi(x,y;w,0) &= x^{-1} \!\!\!\!\!\sum \begin{Sb}d\le x^{1/k} \\ P(d)\le
\min\{y,w\}\end{Sb} \!\!\!\!\!\mu(d)\!\! \sum
\begin{Sb}f\le\sqrt{xd^{-k}} \\ P(f)\le y \\ p(f)> w\end{Sb}
\!\!\mu(f) \sum \begin{Sb}n\le x \\ P(n)\le y \\ d^k\mid n \\ f^2\mid
n\end{Sb} 1 \\
&= \sum \begin{Sb}d\le x^{1/k} \\ P(d)\le \min\{y,w\}\end{Sb}
\!\!\!\!\!\mu(d)\!\! \sum \begin{Sb}f\le\sqrt{xd^{-k}} \\ P(f)\le y \\
p(f)> w\end{Sb} \!\!\mu(f) (x^{-1}\Psi(xd^{-k}\!f^{-2},y)) \\
&= \!\!\!\sum \begin{Sb}d\le x^{1/k} \\ P(d)\le \min\{y,w\}\end{Sb}
\!\!\!\!\mu(d)d^{-k}\!\!\! \sum \begin{Sb}f\le\sqrt{xd^{-k}} \\
P(f)\le y \\ p(f)> w\end{Sb} \!\!\mu(f)f^{-2} \rho \Big(
\lol{xd^{-k}\!f^{-2}}y \Big) \Big({ 1+O\Big( {\lg2x\over\log y} \Big) }\Big),
\end{split}
\end{equation*}
where the final equality follows from Lemma \ref{hild}. The primary
contribution to the double sum will come from those terms with $d$
small and $f=1$ (notice that if $f>1$ then $f>w\ge\log x$), so we
write it as
\begin{multline}
x^{-1}\Psi(x,y;w,0) = \sum_{d<\log x} \!\mu(d)d^{-k} \rho \Big(
\lol{xd^{-k}}y \Big) \Big({ 1 + O\Big( {\lg2x\over\log y} \Big) }\Big) \\
+ O\bigg( \sum_{\log x\le d\le x^{1/k}} \!\!\!d^{-k} \rho \Big(
\lol{xd^{-k}}y \Big) + \sum_{d=1}^\infty d^{-k}\!\!\!\! \sum_{\log x\le
f\le\sqrt x} \!\!\!f^{-2} \rho \Big( \lol{x\!f^{-2}}y \Big) \bigg).
\label{fsboring}
\end{multline}

The two final terms can be estimated by Lemma \ref{rholemma}(c):
\begin{equation*}
\begin{split}
\sum_{\log x\le d\le x^{1/k}} \!\!\!d^{-k} \rho \Big( \lol{xd^{-k}}y \Big)
&\ll \rho\Big( \lol xy \Big) \log^{-k+1}x, \\
\sum_{d=1}^\infty d^{-k}\!\!\!\! \sum_{\log x\le f\le\sqrt x}
\!\!\!f^{-2} \rho \Big( \lol{x\!f^{-2}}y \Big) &\ll \rho\Big( \lol xy \Big)
\log^{-1}x \,\sum_{d=1}^\infty d^{-k} = \zeta(k) \rho\Big( \lol xy
\Big) \log^{-1}x.
\end{split}
\end{equation*}
Also, for $d<\log x$, Lemma \ref{rholemma}(d) gives us
\begin{equation*}
\rho \Big( \lol{xd^{-k}}y \Big) = \rho \Big( \lol xy \Big) \Big({ 1 + O\Big(
{(\log(\log^k\!x))\lgs2^2x\over\log y} \Big) }\Big).
\end{equation*}
Thus equation (\ref{fsboring}) becomes
\begin{equation*}
x^{-1}\Psi(x,y;w,0) = \rho \Big( \lol xy \Big) \bigg( \Big({ 1+O\Big(
\err \Big) }\Big)\! \sum_{d<\log x} \!\mu(d)d^{-k} + O \big( (\log
x)^{-k+1} + \log^{-1}x \big) \bigg).
\end{equation*}
On noting that
\begin{equation*}
\sum_{d<\log x} \mu(d)d^{-k} = \sum_{d=1}^\infty \mu(d)d^{-k} +
O\bigg( \, \sum_{d\ge\log x} d^{-k} \bigg) = \zeta(k)^{-1} +
O(\log^{-k+1}x),
\end{equation*}
part (a) of the lemma is established.

Let $\A_1(x,y;w,0)$ and $\A_2(x,y;w,0)$ denote the 
odd and even elements of $\A(x,y;w,0)$, respectively. For a set $\S$
and an integer $n$, write $n\cdot\S=\{ns: s\in\S\}$ and note that
$\abs{n\cdot\S}=\abs\S$. We see that
\begin{equation*}
\begin{split}
2\cdot\A(x/2,y;w,0) &= \{ 2n: n\le x/2;\, n\hbox{ is $k$-free};\,
P(n)\le y;\, d^2\mid n\implies P(d)\le w \} \\
&= \{ n: n\le x;\, n\hbox{ is even and $k$-free};\, P(n)\le y;\,
d^2\mid n\implies P(d)\le w \} \\
&\quad\,\cup \{ 2^kn: n\le x/2^k;\, n\hbox{ is odd and $k$-free};\,
P(n)\le y;\, d^2\mid n\implies P(d)\le w \} \\
&= \A_2(x,y;w,0) \cup 2^k\cdot\A_1(x/2^k,y;w,0)
\end{split}
\end{equation*}
as a disjoint union. Taking cardinalities on both sides, we see that
\begin{equation*}
\Psi(x/2,y;w,0) = \Psi_2(x,y;w,0) + \Psi_1(x/2^k,y;w,0),
\end{equation*}
or equivalently
\begin{equation*}
\Psi_1(x,y;w,0) - \Psi_1(x/2^k,y;w,0) = \Psi(x,y;w,0) -
\Psi(x/2,y;w,0),
\end{equation*}
since $\Psi=\Psi_1+\Psi_2$. Consequently, for any nonnegative integer
$m$, we have
\begin{equation}
\begin{split}
\Psi_1(x,y;w,0) - \Psi_1(x/2^{(m+1)k},y;w,0) &= \sum_{j=0}^m \big(
\Psi_1(x/2^{jk},y;w,0) - \Psi_1(x/2^{(j+1)k},y;w,0) \big) \\
&= \sum_{j=0}^m \big( \Psi(x/2^{jk},y;w,0) - \Psi(x/2^{jk+1},y;w,0)
\big).
\end{split}
\label{refiter}
\end{equation}

Choose $m=\floor{\lg2x/(k\log 2)}$, so that $2^{mk}\le\log
x<2^{(m+1)k}$. For any $0\le l\le(m+1)k$, part (a) gives us
\begin{equation*}
\begin{split}
\Psi(x/2^l,y;w,0) &= \frac x{2^l\zeta(k)} \rho \Big( \lol{(x/2^l)}y
\Big) \Big( 1 + O\Big( \err \Big) \Big) \\
&= \frac x{2^l\zeta(k)} \rho \Big( \lol xy \Big) \Big( 1 + O\Big(
{\log(2^k\log x)\lgs2^2x \over \log y} \Big) \Big) \Big( 1 + O\Big(
\err \Big) \Big)
\end{split}
\end{equation*}
by Lemma \ref{rholemma}(d) and the choice of $m$. Using this in
equation (\ref{refiter}) gives us
\begin{equation*}
\begin{split}
\Psi_1(x,y;w,0) &= \sum_{j=0}^m \Big( \frac x{2^{jk}\zeta(k)} - \frac
x{2^{jk+1}\zeta(k)} \Big) \rho \Big( \lol xy \Big) \Big( 1 + O\Big(
\err \Big) \Big) \\
&\qquad+ O(\Psi(x/2^{(m+1)k},y;w,0)) \\
&= \frac x{2\zeta(k)} \rho \Big( \lol xy \Big) \Big( 1 + O\Big(
\err \Big) \Big) \sum_{j=0}^m \frac1{2^{jk}} + O\Big( \frac
x{2^{(m+1)k}} \rho \Big( \lol xy \Big) \Big) \\
&= \frac x{2\zeta(k)} \rho \Big( \lol xy \Big) \Big( 1 + O\Big(
\err \Big) \Big) \big( (1-2^{-k})^{-1} + O\big( \frac1{\log x} \big)
\big) \\
&\qquad+ O\Big( \frac x{\log x} \rho \Big( \lol xy \Big) \Big) \\
&= \frac x{\xi(k)} \rho \Big( \lol xy \Big) \Big( 1 + O\Big(
\err \Big) \Big).
\end{split}
\end{equation*}

Finally, the number of integers not exceeding $x$ of the form
$m^2+m-1$ is at most $\sqrt x$. Therefore,
\begin{equation*}
\begin{split}
\Psi_0(x,y;w,0) &= \Psi_1(x,y;w,0) + O(\sqrt x) \\
&= \frac x{\xi(k)} \rho \Big( \lol xy \Big) \Big( 1 + O\Big( \err +
x^{-1/2}\rho\Big( \lol xy \Big)^{\!\!-1} \Big) \Big).
\end{split}
\end{equation*}
Since the second error term is $\ll x^{-1/2+\ep}$ by Lemma
\ref{rholemma}(b), it is dominated by the first error term, and so
part (b) is established.
\end{proof}

\begin{lemma}
Let $k\ge2$ be an integer and $0<\ep<1/k$ and $0<\lambda<1$ be real
numbers. Let $x$ be a sufficiently large real number, and let $y$ and
$w$ be real numbers and $p$ a prime satisfying $\log x\le w\le
x^{(1-\ep)/k}$ and $\exp(\lgs2^4x)\le p<y\le x^{(1-\ep)/2}$. Let $l<k$ a
positive integer, with $l=1$ if $p\ge w$. Then:
\begin{enumerate}
\item[\rm(a)]{
$\displaystyle \Psi(x,y;w,\lambda) = (1-\lambda)\frac x{\zeta(k)} \rho
\Big( \lol xy \Big) \Big({ 1 + O\Big( \err \Big) }\Big)$;
}
\item[\rm(b)]{
$\displaystyle \Psi_0(x,y;w,\lambda) = (1-\lambda)\frac x{\xi(k)} \rho
\Big( \lol xy \Big) \Big({ 1 + O\Big( \err \Big) }\Big)$;
}
\item[\rm(c)]{
$\Psi(x,y;w,\lambda;p^l)>p$, and if $p\ne2$, then
$\Psi_0(x,y;w,\lambda;p^l)>p$;
}
\item[\rm(d)]{
there exists an element of $\A(x,y;w,\lambda)$ that is exactly
divisible by $2^l$ (for any $1\le l<k$).
}
\end{enumerate}
\label{smoothsets}
\end{lemma}

\begin{proof}
Clearly $\Psi(x,y;w,\lambda)=\Psi(x,y;w,0)-\Psi(\lambda x,y;w,0)$ and
similarly for the $\Psi_0$, and so parts (a) and (b) follow from Lemma
\ref{kfrees} on noting that
\begin{equation*}
\rho \Big( \lol{\lambda x}y \Big) = \rho \Big( \lol xy \Big) \Big({ 1
+ O\Big( \err \Big) }\Big)
\end{equation*}
by Lemma \ref{rholemma}(d) with $v=\lambda^{-1}$.

We notice that $mp^l\in\A(x,y;w,\lambda;p^l)$ if and only if
$m\in\A(xp^{-l},p-1;w,\lambda)$, and so
\begin{equation*}
\Psi(x,y;w,\lambda;p^l) = \Psi(xp^{-l},p-1;w,\lambda) = \rho \Big(
\lol{xp^{-l}}{(p-1)} \Big) {(1-\lambda)xp^{-l} \over \zeta(k)} \Big({
1 + O\Big( \frac{\lgs2^3x}{\log(p-1)} \Big) }\Big)
\end{equation*}
by part (a). To show that this is greater than $p$, it suffices to
show that
\begin{equation}
{\zeta(k)p^{l+1}\over1-\lambda} \rho \Big( \lol{xp^{-l}}{(p-1)}
\Big)^{\!\!-1} < x \Big({ 1 + O\Big( \frac{\lgs2^3x}{\log(p-1)} \Big) }\Big).
\label{collectps}
\end{equation}
But under the restrictions on $p$, $l$, $y$, and $w$, we have that
$p^{l+1}<x^{1-\ep}$; and by Lemma \ref{rholemma}(b) we have that $\rho
( \lol{xp^{-l}}{(p-1)} )^{-1} \le \rho ( \lol x{(p-1)} )^{-1} \ll
x^{\ep/2}$. Therefore, the left-hand side of equation
(\ref{collectps}) is $\ll x^{1-\ep/2}$. This establishes the first
assertion of part (c), and the proof of the second assertion is
similar.

Finally, by part (b) with $x$ replaced by $x/2^l$, we see that
$\Psi_0(x/2^l,y;w,\lambda)$ is positive for sufficiently large
$x$. But then we can choose an element $m$ of
$\A_0(x/2^l,y;w,\lambda)$, and then $2^lm$ is an element of
$\A(x,y;w,\lambda)$ that is exactly divisible by $2^l$, which
establishes part (d).
\end{proof}

\begin{lemma}
Under the hypotheses of Lemma \ref{smoothsets}, we have
\begin{enumerate}
\item[\rm(a)]{$\displaystyle
\sum_{n\in\A(x,y;w,\lambda)} \!\!\!\!n^{-1} = \rho \Big( \lol xy \Big)
\frac{\log\lambda^{-1}}{\zeta(k)} \Big({ 1 + O\Big( \err \Big) }\Big)
$;}
\item[\rm(b)]{$\displaystyle
\sum_{n\in\A_0(x,y;w,\lambda)} \!\!\!\!n^{-1} = \rho \Big( \lol xy \Big)
\frac{\log\lambda^{-1}}{\xi(k)} \Big({ 1 + O\Big( \err \Big) }\Big)
$;}
\item[\rm(c)]{if $y=y(x)$ is chosen so that $\log x/\log y\le B$ for
some constant $B$ depending only on $k$, then for any positive real
number $\alpha$, there exists a real number $0<\lambda<1$, bounded
away from zero uniformly in $x$, such that when $x$ is sufficiently
large, we have
\begin{equation}
0<\alpha-\!\!\!\!\sum_{n\in\A_0(x,y;w,\lambda)}\!\!\!\!n^{-1}<x^{-1}
\exp(D\alpha)
\label{explicit}
\end{equation}
for some constant $D$ depending only on $k$.}
\end{enumerate}
\label{recipsums}
\end{lemma}

\vskip-8pt\noindent When we apply part (c) we will need only that the
middle expression in equation (\ref{explicit}) is positive and $\ll
x^{-1}$, but we have been explicit about the upper bound for two
reasons. First, since we will use part (c) to choose a value of
$\lambda$, it is important that the implicit constant not depend on
$\lambda$. Second, we will apply part (c) with a very small value of
$\alpha$, and one that is specified only up to an error bounded by a
function of $x$. From the form of the upper bound in equation
(\ref{explicit}), we will be assured that this will not inflate the
implicit constant as $x$ grows.\vskip12pt

\begin{proof}
We have
\begin{equation}
\begin{split}
\zeta(k) \!\!\!\!\sum_{n\in\A(x,y;w,\lambda)} \!\!\!\!n^{-1} &= \zeta(k)
\int_{\lambda x}^x t^{-1} \,d\Psi(t,y;w,0) \\
&= \zeta(k)\, t^{-1}\Psi(t,y;w,0)\big\vert_{\lambda x}^x + \zeta(k)
\int_{\lambda x}^x \Psi(t,y;w,0) \,{dt\over t^2} \\
&= \bigg({ \rho \Big( \lol xy \Big) - \rho \Big( \lol{\lambda x}y
\Big) + \int_{\lambda x}^x \rho \Big( \lol ty \Big) \frac{dt}t }\bigg)
\Big({ 1 + O\Big( \err \Big) }\Big)
\end{split}
\label{lxtox}
\end{equation}
by Lemma \ref{kfrees}(a). However, since $\log\lambda^{-1}=O(1)=o(\log
y/\lgs2^2x)$, we can apply Lemma \ref{rholemma}(d) to see that
\begin{equation*}
\rho \Big( \lol ty \Big) = \rho \Big( \lol xy \Big) \Big({ 1 + O\Big(
\err \Big) }\Big)
\end{equation*}
uniformly for $\lambda x\le t\le x$. Therefore equation (\ref{lxtox})
becomes
\begin{equation*}
\begin{split}
\zeta(k) \!\!\!\!\sum_{n\in\A(x,y;w,\lambda)} \!\!\!\!n^{-1} &= \rho
\Big( \lol xy \Big) \bigg({ O\Big( \err \Big) + \int_{\lambda x}^x
\frac{dt}t }\bigg) \Big({ 1 + O\Big( \err \Big) }\Big) \\
&= \rho \Big( \lol xy \Big) \log \lambda^{-1} \Big({ 1 + O\Big( \err
\Big) }\Big).
\end{split}
\end{equation*}
This establishes part (a), and the proof of part (b) is exactly
similar.

As for part (c), note that the function
$\alpha-\sum_{n\in\A_0(x,y;w,t)}n^{-1}$ is an increasing function of
$t$, with jump discontinuities of size not exceeding $(tx)^{-1}$, and
which takes negative values if $x$ is sufficiently large. Therefore we
can choose $\lambda$ such that
\begin{equation}
0<\alpha-\!\!\!\!\sum_{n\in\A_0(x,y;w,\lambda)}\!\!\!\!n^{-1}
\le(\lambda x)^{-1}.
\label{lelaminv}
\end{equation}
For this value of $\lambda$, part (b) tells us that
\begin{equation*}
\alpha>\!\!\!\!\sum_{n\in\A_0(x,y;w,\lambda)}\!\!\!\!n^{-1} =
{\log\lambda^{-1}\over\xi(k)} \rho \Big( \lol xy \Big) \Big({ 1 +
O\Big( \err \Big) }\Big),
\end{equation*}
which implies that
\begin{equation*}
\lambda^{-1} < \exp\bigg({ \alpha {\xi(k)\over\rho(B)} \Big({ 1 +
O\Big( \err \Big) }\Big) }\bigg).
\end{equation*}
This shows that $\lambda$ is bounded away from zero uniformly in $x$;
and combining this bound with equation (\ref{lelaminv}) and writing
$D=2\xi(k)/\rho(B)$ establishes equation (\ref{explicit}) for
sufficiently large $x$.
\end{proof}

\section{Construction of Dense Egyptian Fractions}

We now proceed with the proof of Theorem \ref{densethm}. Let $r=a/b$
be a positive rational number and $\eta$ a positive real number. For
the first stage of the proof, let $0<\ep<\half$ and $0<\delta<r$ be
real numbers and $k\ge2$ an integer that is large enough to ensure
that $b$ is $k$-free; we will later constrain these parameters in
terms of $r$ and $\eta$. Set
\begin{equation}
\lambda = \exp\Big( -{(r-\delta)\zeta(k)\over\rho(2/(1-\ep))} \Big).
\label{lamdef}
\end{equation}
Let $x\ge3$ be a sufficiently large real number, set $y=x^{(1-\ep)/2}$
and $w=x^{(1-\ep)/k}$, and for $z>2$ define
\begin{equation*}
\begin{split}
D(z) &= D(z;w) = \bigg( \prod \begin{Sb}p<z \\ p\le
w\end{Sb}\!p^{k-1} \bigg) \bigg( \prod_{w<p<z}\!p \bigg), \\
D_0(z) &= 2^{-k+1}D(z),
\end{split}
\end{equation*}
so that $D_0(z)$ is the odd part of $D(z)$. (We understand that the
second product in the definition of $D(z)$ is 1 if $z\le w$.) Notice
that if $p$ is a prime, then $p$ does not divide $D(p)$.

Let $p_1>p_2>\dots>p_R$ be the primes in $(w,y]$, and for later
consistency of notation, let $p_0$ be the smallest prime exceeding
$y$. Notice that all elements of $\A(x,y;w,0)$ divide $D(p_0)$, as
does $b$ as long as $P(b)\le w$; we henceforth assume that $x$ is
large enough to ensure this. Define
\begin{equation*}
\abfrac0 =  r -\!\!\! \sum_{n\in\A(x,y;w,\lambda)} \!\frac1n
\end{equation*}
with $\abfrac0$ in lowest terms. Notice that $b_0$ divides $D(p_0)$
and that
\begin{equation}
\abfrac0 = \delta + O\Big( \err \Big) \label{arbrsize}
\end{equation}
by Lemma \ref{recipsums}(a) and the choice (\ref{lamdef}) of $\lambda$.

We now recursively construct a sequence of fractions $\abfrac1$,
$\abfrac2$, \dots, $\abfrac R$ with the following two
properties. First, each $\abfrac i$ is obtained from the previous
$\abfrac{i-1}$ by adding the reciprocals of a few elements of
$\A(x,y;w,\lambda)$, specifically elements of $\A(x,y;w,\lambda;p_i)$;
so when we have constructed these fractions, we will have written $r$
as the sum of $\abfrac R$ and the reciprocals of almost all of the
elements of $\A(x,y;w,\lambda)$. Second, all of the primes dividing
each denominator $b_i$ will be less than $p_i$, so the denominators
are becoming gradually smoother.

Formally, we construct the fractions $\abfrac i$ as follows.  Given
$\abfrac{i-1}$, where $1\le i\le R$, we apply Lemma \ref{addfracs}
with $p=p_i$, $l=1$, $N=D(p_{i-1})$, $c/d=\abfrac{i-1}$, and
$\S=\A(x,y;w,\lambda;p_i)$.  This gives us a subset of
$\A(x,y;w,\lambda;p_i)$, which we call $\B_i$, with cardinality less
than $p_i$ such that, if we define
\begin{equation*}
\abfrac i = \abfrac{i-1} + \sum_{n\in\B_i} \frac1n
\end{equation*}
with $\abfrac i$ in lowest terms, then $b_i$ divides
$D(p_{i-1})/p_i=D(p_i)$. We are justified in applying Lemma
\ref{addfracs} with these parameters since, by Lemma
\ref{smoothsets}(c), the size of $\A(x,y;w,\lambda;p_i)$ is at
least~$p_i$.

Now set $y'=\exp(\lgs2^4 x)$, and let $q_1>q_2>\dots>q_S$ be the
primes in $[y',w]$; and for later consistency of notation, set
$q_0=p_R$ and $\abcfrac0k=\abfrac R$. We recursively construct another
sequence of fractions $\abcfrac11$, $\abcfrac12$, \dots, $\abcfrac1k$,
$\abcfrac21$, \dots, $\abcfrac Sk$ with similar properties to the
first sequence, the differences coming from the fact that the primes
$q_i$ may potentially divide the denominators we are working with to the
$(k-1)$st power. First, each $\abcfrac i1$ is equal to
$\abcfrac{i-1}k$, and each $\abcfrac ij$ is obtained from the previous
$\abcfrac i{j-1}$ by adding the reciprocals of a few elements of
$\A(x,y;w,\lambda)$, specifically elements of
$\A(x,y;w,\lambda;q_i^{k-j+1})$; again, we will have written $r$ as be
the sum of $\abcfrac Sk$ and the reciprocals of almost all of the
elements of $\A(x,y;w,\lambda)$. Second, all of the primes dividing
each denominator $b_{i,j}$ will be less than or equal to $q_i$, and
$q_i$ itself will divide $b_{i,j}$ at most to the $(k-j+1)$st power,
so the denominators are becoming gradually smoother still.

Formally, we construct the fractions $\abcfrac ij$ as follows. First,
given $\abcfrac{i-1}k$, where $1\le i\le S$, we set $\abcfrac
i1=\abcfrac{i-1}k$. Then, given $\abcfrac i{j-1}$, where $2\le j\le
k$, we apply Lemma \ref{addfracs} with $p=q_i$, $l=k-j+1$,
$N=q_i^{k-j+1}D(q_i)$, $c/d=\abcfrac i{j-1}$, and
$\S=\A(x,y;w,\lambda;q_i^{k-j+1})$. This gives us a subset of
$\A(x,y;w,\lambda;q_i^{k-j+1})$, which we call $\B_{i,j}$, with
cardinality less than $q_i$ such that, if we define
\begin{equation*}
\abcfrac ij = \abcfrac i{j-1} + \sum_{n\in\B_{i,j}} \frac1n
\end{equation*}
with $\abcfrac ij$ in lowest terms, then $b_{i,j}$ divides
$q_i^{k-j}D(q_i)$. Again, Lemma \ref{smoothsets}(c) justifies our use
of Lemma \ref{addfracs} with these parameters. We also note that
$D(q_{i-1})=q_i^{k-1}D(q_i)$, and so our convention that
$\abcfrac{i-1}k=\abcfrac i1$ is consistent with the divisibility
property of the denominators.

We have gradually eliminated all prime factors not less than $y'$ from
the denominators of the fractions $\abcfrac ij$; to finish the first
stage of the proof, we eliminate powers of two as well. We set
$\abcfrac{S+1}1=\abcfrac Sk$ and define a sequence of fractions
$\abcfrac{S+1}2$, $\abcfrac{S+1}3$, \dots, $\abcfrac{S+1}k$, such that
each $\abcfrac{S+1}j-\abcfrac{S+1}{j-1}$ is either zero or else is the
reciprocal of an element of $\A(x,y';w,\lambda)$, and such that powers
of two dividing the $b_{S+1,j}$ are diminishing. Given
$\abcfrac{S+1}{j-1}$, where $2\le j\le k$, if $2^{k-j+1}$ exactly
divides $b_{S+1,j-1}$ then we invoke Lemma \ref{smoothsets}(d) to
choose an element $n_j\in\A(x,y';w,\lambda)$ that is exactly divisible
by $2^{k-j+1}$ and set $\abcfrac{S+1}j=\abcfrac{S+1}{j-1}+1/n_j$ and
$\B_{S+1,j}=\{n_j\}$. Otherwise, we set
$\abcfrac{S+1}j=\abcfrac{S+1}{j-1}$ and $\B_{S+1,j}=\emptyset$. It is
easy to see inductively that with this construction, $b_{S+1,j}$
divides $2^{k-j}D_0(y')$ for each $1\le j\le k$.

Define
\begin{equation}
\B=\bigg( \bigcup_{i=1}^R \B_i \bigg) \cup \bigg( \bigcup_{i=1}^{S+1}
\bigcup_{j=2}^k \B_{i,j} \bigg).
\label{Bunion}
\end{equation}
We note the largest prime factor of each element of a $\B_i$ is $p_i$,
and the largest prime factor of each element of a $\B_{i,j}$ is $q_i$,
except that the largest prime factor of each element of a $\B_{S+1,j}$
is less than $y'$. Furthermore, for a fixed $i\le S$, the various
$\B_{i,j}$ have elements which are divisible by $q_i$ to different
powers; and the elements of the various $\B_{S+1,j}$ are divisible by
different powers of two. We conclude that the union in the definition
(\ref{Bunion}) of $\B$ is actually a disjoint union, and thus if
we define $\A=\A(x,y;w,\lambda)\setminus\B$, then we have written
\begin{equation}
r=\abcfrac{S+1}k+\sum_{n\in\A}\frac1n \label{firststage}
\end{equation}
with $b_{S+1,k}$ dividing $D_0(y')$. Moreover, the cardinality of $\B$
is bounded by
\begin{equation*}
\sum_{w<p\le y}(p-1) + (k-1)\sum_{y'\le p\le w}(p-1) + (k-1) \ll y^2 =
x^{1-\ep},
\end{equation*}
so that the cardinality of $\A$ is
\begin{equation}
\begin{split}
\abs{\A} &\ge \abs{\A(x,y;w,\lambda)} - O(x^{1-\ep}) \\
&= (1-\lambda)\frac x{\zeta(k)} \rho\Big( {2\over1-\ep} \Big) \Big({ 1 +
O\Big( \err \Big) }\Big)
\end{split}
\label{carda}
\end{equation}
by Lemma \ref{smoothsets}(a). We also note that
\begin{equation}
\abcfrac{S+1}k = \abfrac0+\sum_{n\in\B} \frac1n = \Big({ \delta +
O\Big( \err \Big) }\Big) + O\big( \abs\B\!(\lambda x)^{-1} \big) =
\delta + O\Big( \err \Big)
\label{notfine}
\end{equation}
by equation (\ref{arbrsize}) and the fact that $\abs\B\ll x^{1-\ep}$.

From the choice (\ref{lamdef}) of $\lambda=\lambda(\ep,\delta,k)$ and
the continuity of $\rho(u)$, we see that
\begin{equation*}
\lim\begin{Sb}\ep\to0 \\ \delta\to0 \\ k\to\infty\end{Sb}
{1-\lambda\over\zeta(k)} \rho\Big( {2\over1-\ep} \Big) = \Big({
1-\exp\Big( \frac {-r}{\rho(2)} \Big) }\Big) \rho(2) = C(r).
\end{equation*}
We thus choose $\ep$ and $\delta$ sufficiently small and $k$
sufficiently large, in terms of $r$ and $\eta$, to ensure that
$(1-\lambda)\rho(2/(1-\ep))/\zeta(k) > C(r)-\eta$; then, from equation
(\ref{carda}), we see that $\abs\A>(C(r)-\eta)x$ for sufficiently
large $x$. (We note that now $x$ needs only to be sufficiently large
in terms of $r$ and $\eta$, since $\delta$, $\ep$, and $k$, and thus
$\lambda$, have all been chosen in terms of $r$ and $\eta$.) Moreover,
all elements of $\A$ are certainly greater than $\lambda
x$. Therefore, to establish Theorem~\ref{densethm}, it suffices by
equation (\ref{firststage}) to write $\abcfrac{S+1}k$ as the sum of
reciprocals of distinct integers not exceeding $\lambda x$, without
regard to the number of terms in the representation. This is the goal
of the second stage of the proof.

We begin by applying to $\abcfrac{S+1}k$ much the same process that we
applied to $r$ in the first stage of the proof. Recall that $y'=\exp(\lgs2^4x)$, and set $x'=(y')^{2k}$, $w'=y'$, and $y''=(\log
x)/3k$. Let $q'_1>q'_2>\dots>q'_T$ be the primes in $[y'',y')$, and
for later consistency of notation, set $q'_0=q_S$. Since $\log x'/\log
y'=2k$ and $\exp(\lgs2^4x')\le y''<y'=w'<(x')^{(1-\ep)/k}$ for $x$
sufficiently large, we may freely appeal to Lemmas \ref{smoothsets}
and \ref{recipsums} with the parameters $x'$, $y'$, and $w'$ and any
prime $q'_i$. In particular, by Lemma \ref{recipsums}(c), we may
choose $\lambda'$ so that
\begin{equation}
0<\abcfrac{S+1}k-\!\!\!\sum_{n\in\A_0(x',y';y',\lambda')}\!\!\!\!\!n^{-1}
< (x')^{-1} \exp\Big( D(k)\abcfrac{S+1}k \Big).
\label{buttpain}
\end{equation}

Define
\begin{equation*}
\cdfrac0k=\abcfrac{S+1}k-\!\!\sum_{n\in \A_0(x',y'; y',\lambda')}\frac1n
\end{equation*}
with $\cdfrac0k$ in lowest terms, and notice that $b'_{0,k}$
divides $D_0(y')$. With this definition and the estimate
(\ref{notfine}) on the size of $\abcfrac{S+1}k$, equation
(\ref{buttpain}) becomes
\begin{equation*}
0<\cdfrac0k<(x')^{-1}\exp\bigg({ D(k) \Big({ \delta+O\Big( \err \Big)
}\Big) }\bigg) \ll (x')^{-1}.
\end{equation*}

We recursively construct a sequence of fractions $\cdfrac11$,
$\cdfrac12$, \dots, $\cdfrac1k$, $\cdfrac21$, \dots, $\cdfrac Tk$ as
follows. First, given $\cdfrac{i-1}k$, where $1\le i\le T$, we set
$\cdfrac i1=\cdfrac{i-1}k$. Then, given $\cdfrac i{j-1}$, where $2\le
j\le k$, we apply Lemma \ref{addfracs} with $p=q'_i$, $l=k-j+1$,
$N=(q'_i)^{k-j+1}D(q'_i)$, $c/d=\cdfrac i{j-1}$, and
$\S=\A(x',y';y',\lambda';(q'_i)^{k-j+1})$. This gives us a subset of
$\A(x',y';y',\lambda';(q'_i)^{k-j+1})$, which we call $\B'_{i,j}$, with
cardinality less than $q'_i$ such that, if we define
\begin{equation*}
\cdfrac ij = \cdfrac i{j-1} + \sum_{n\in\B'_{i,j}} \frac1n
\end{equation*}
with $\cdfrac ij$ in lowest terms, then $b'_{i,j}$ divides
$(q'_i)^{k-j}D(q'_i)$.

Define
\begin{equation*}
\B' = \bigcup_{i=1}^T \bigcup_{j=2}^k \B'_{i,j},
\end{equation*}
and notice that this is a disjoint union by the reasoning following
equation (\ref{Bunion}) earlier. Thus, if we define
$\A'=\A_0(x',y';y',\lambda')\setminus\B'$, we see that
\begin{equation}
\abcfrac{S+1}k = \cdfrac Tk+ \sum_{n\in\A'} \frac1n
\end{equation}
with $b'_{T,k}$ dividing $D_0(y'')$; in particular,
$P(b'_{T,k})<y''$. Moreover,
\begin{equation*}
0<\cdfrac Tk = \cdfrac0k + \sum_{n\in\A'} \frac1n \ll (x')^{-1} +
\abs{\A'}\!(\lambda'x')^{-1} \ll (y')^2(x')^{-1} = o((y'')^{-1}),
\end{equation*} 
and so $\cdfrac Tk<(P(b'_{T,k}))^{-1}$ for $x$ sufficiently
large. Since $b'_{T,k}$ divides $D_0(y'')$, it is odd, and so we can 
apply Lemma \ref{Breusch} to find a set $\C$ of positive odd integers
such that $\cdfrac Tk= \sum_{n\in\C} 1/n$ and
\begin{equation*}
\max\{n\in\C\} \ll b'_{T,k} \!\!\!\prod_{p\le P(b'_{T,k})}\!\! p \le
\prod_{p<y''} p^k = \exp\!\bigg( k\sum_{p<y''}\log p \bigg).
\end{equation*}

Chebyshev's bound for $\pi(t)$ implies that there is a real number
$c<2$ such that $\sum_{p<t}\log p<ct$ for any positive
$t$. Therefore
\begin{equation*}
\max\{n\in\C\} \ll \exp(kcy'') < x^{2/3}.
\end{equation*}

We have almost achieved our goal for the second stage of the proof,
for $\abcfrac{S+1}k=\sum_{n\in\A'} 1/n + \sum_{n\in\C} 1/n$; but $\A'$
and $\C$ might not be disjoint. Define
\begin{equation*}
\D_1 = \{n+1: n\in\A'\cap\C\} \quad\hbox{and}\quad \D_2 = \{n(n+1):
n\in\A'\cap\C\}.
\end{equation*}
Since $1/n=1/(n+1)+1/(n(n+1))$, we have
\begin{equation}
\begin{split}
\abcfrac{S+1}k &= \sum_{n\in\A'}\frac1n +
\sum_{n\in\C\setminus\A'}\frac1n + \sum_{n\in\A'\cap\C}\frac1n \\
&= \sum_{n\in\A'}\frac1n + \sum_{n\in\C\setminus\A'}\frac1n +
\sum_{n\in\D_1}\frac1n + \sum_{n\in\D_2}\frac1n,
\end{split}
\label{foursets}
\end{equation}
and we claim that this is a disjoint representation of
$\abcfrac{S+1}k$ using denominators less than $\lambda x$. We already
know that $\max\{n\in\A'\} \le x'$ and $\max\{n\in\C\} \ll x^{2/3}$,
and we easily see that
\begin{equation*}
\max\{n\in\D_1\cup\D_2\} \ll (\max\{n\in\A'\})^2 \le (x')^2,
\end{equation*}
so that the integers involved are of admissible size.

Clearly $\A'$ and $\C\setminus\A'$ are disjoint; and since the
elements of $\A'$ and $\C$ are odd, those of $\D_1$ and $\D_2$ are
even, and so each of the first two is disjoint from each of the last
two. Finally, if there were an element $n$ in $\D_1\cap\D_2$, then
there would exist $m_1$, $m_2\in\A'\cap\C$ such that
$n=m_1+1=m_2(m_2+1)$. But then
$m_1\in\A'\subset\A_0(x',y';y',\lambda')$ satisfies $m_1=m_2^2+m_2-1$,
contradicting the definition of an $\A_0$-set. Therefore the four sets
in the representation (\ref{foursets}) are indeed disjoint, and
Theorem \ref{densethm} is established.

\section{Prospects for Improved Results}

Clearly the only barrier to establishing Theorem \ref{densethm} in
best-possible form is the presence of the factors $1-\log2$ in the
expression for $C(r)$. This quantity arises as $\rho(2)$, which comes
from using integers of size $x$ that are roughly $x^{1/2}$-smooth,
which in turn comes from the necessity that we have at least $p-1$
multiples of every prime $p$ at our disposal in order to invoke Lemma
\ref{Zpsums} and thus Lemma \ref{addfracs}.

As noted in Section \ref{elemsec}, the conclusion of Lemma
\ref{Zpsums} is best possible, since we can arrange for many sums of
subsets to coincide.  For randomly chosen sets, however, this behavior
is very unlikely. Heuristically, a set of random nonzero elements of
${\bf Z}_p$ should have a subset summing to a randomly chosen target
element of ${\bf Z}_p$ as soon as the size of the set is as large as a
small power of $\log p$. We might expect, then, that if we were to
attempt the construction in the proof of Theorem \ref{densethm} using
integers that were $x^{1-\ep}$-smooth instead of
$x^{(1-\ep)/2}$-smooth, then we would have at least $p^\ep$ multiples
of each prime $p$ to choose from, and almost always we could find a
small subset of those multiples to exclude to force the necessary
cancellation of factors of $p$. If one could show that this were the
case, the factors tending to $\rho(2)$ could be replaced by factors
tending to $\rho(1)=1$, and we would have established the best
possible result.

Number theorists have also investigated whether positive rationals
have Egyptian fraction representations of various specified forms. For
instance, we have seen in Lemma \ref{Breusch} that any positive
rational with odd denominator can be written as the sum of reciprocals
of distinct odd positive integers (clearly no such representation
exists for positive rationals whose denominator in reduced form is
even, since the common denominator of such a representation will
necessarily be odd). Graham \cite{Gra:OFSoUF} has proven a very
general theorem showing that, for a certain class of subsets $\Z$ of
the positive integers, any positive rational can be written as an
Egyptian fraction with denominators restricted to elements of $\Z$,
provided only that it satisfies to clearly necessary assumptions. Its
denominator must not be divisible by any primes that don't divide any
element of $\Z$, as in the case of odd denominators discussed above;
and its size must be compatible with the sizes of the finite subsums
of the reciprocals of $\Z$---the series of reciprocals of $\Z$ might
converge, for instance, in which case large rationals could never be
so represented.

It does not seem implausible, therefore, that any positive rational
with odd denominator has an Egyptian fraction representation
consisting of $\gg x$ unit fractions all of whose denominators are odd
and at most $x$; or that a positive rational meeting the local
conditions prescribed by $\Z$ can be written as an Egyptian fraction
using an asymptotically positive proportion of the elements of
$\Z$. Furthermore, one might even believe that the best-possible bound
for that proportion, which can be derived from the function
$\sum_{n\in\Z,\,\lambda x<n\le x}1/n$ as we did in equation
(\ref{derive}) for $\Z={\bf Z}$, is in fact attainable.

\providecommand{\bysame}{\leavevmode\hbox to3em{\hrulefill}\thinspace}

\end{document}